\documentclass{amsart}

\newtheorem{thm*}{Theorem}

\newtheorem{thm}{Theorem}[section]

\newtheorem{lemma}{Lemma}

\newtheorem{remark}{Remark}

\newtheorem{cor}{Corollary}

\numberwithin{equation}{section}

\begin{document}

\def\L{ {\mathcal{L}}}
\def\d{ \partial }
\def\Na{{\mathbb{N}}}

\def\Z{{\mathbb{Z}}}

\def\IR{{\mathbb{R}}}

\newcommand{\E}[0]{ \varepsilon}

\newcommand{\la}[0]{ \lambda}

\newcommand{\s}[0]{ \mathcal{S}}

\newcommand{\AO}[1]{\| #1 \| }

\newcommand{\BO}[2]{ \left( #1 , #2 \right) }

\newcommand{\CO}[2]{ \left\langle #1 , #2 \right\rangle}

\newcommand{\R}[0]{ \IR\cup \{\infty \} }

\newcommand{\co}[1]{ #1^{\prime}}

\newcommand{\p}[0]{ p^{\prime}}

\newcommand{\m}[1]{   \mathcal{ #1 }}

%  norm of H
\newcommand{ \A}[1]{ \left\| #1 \right\|_H }

% inner product H
\newcommand{\B}[2]{ \left( #1 , #2 \right)_H }

% H^* , H pairing
\newcommand{\C}[2]{ \left\langle #1 , #2 \right\rangle_{  H^* , H } }

 \newcommand{\HON}[1]{ \| #1 \|_{ H^1} }

% Omega    \Om
\newcommand{ \Om }{ \Omega}

% \partial Omega      \pOm
\newcommand{ \pOm}{\partial \Omega}

%   D(\Omega)   \D
\newcommand{\D}{ \mathcal{D} \left( \Omega \right)}

% D'( Omega)        \DP
\newcommand{\DP}{ \mathcal{D}^{\prime} \left( \Omega \right)  }

% D' pairing
\newcommand{\DPP}[2]{   \left\langle #1 , #2 \right\rangle_{  \mathcal{D}^{\prime}, \mathcal{D} }}

% (H^1)^* , H^1    (( pairing ))      \PHH
\newcommand{\PHH}[2]{    \left\langle #1 , #2 \right\rangle_{    \left(H^1 \right)^*  ,  H^1   }    }

%  H^{-1} , H_0^1  (( pairing ))   \PHO
\newcommand{\PHO}[2]{  \left\langle #1 , #2 \right\rangle_{  H^{-1}  , H_0^1  }}

 %  H^1(\Omega)     \HO
 \newcommand{\HO}{ H^1 \left( \Omega \right)}

%  H_0^1( \Omega)       \HOO
\newcommand{\HOO}{ H_0^1 \left( \Omega \right) }

% C_c^\infty(omega)
\newcommand{\CC}{C_c^\infty\left(\Omega \right) }

%H_0^1(Omega)  norm
\newcommand{\N}[1]{ \left\| #1\right\|_{ H_0^1  }  }

%H_0^1(Omega)   innerproduct
\newcommand{\IN}[2]{ \left(#1,#2\right)_{  H_0^1} }

% H^1(\Omega) inner product
\newcommand{\INI}[2]{ \left( #1 ,#2 \right)_ { H^1}}

% (H^1(\Omega))^*
\newcommand{\HH}{   H^1 \left( \Omega \right)^* }

% ( H^{-1}(\Omega))
\newcommand{\HL}{ H^{-1} \left( \Omega \right) }

\newcommand{\HS}[1]{ \| #1 \|_{H^*}}

\newcommand{\HSI}[2]{ \left( #1 , #2 \right)_{ H^*}}

\newcommand{\Ov}{ \overline{ \Omega}}
\newcommand{\WO}{ W_0^{1,p}}
\newcommand{\w}[1]{ \| #1 \|_{W_0^{1,p}}}

\newcommand{\ww}{(W_0^{1,p})^*}

%%%%%%%%%%%%%%%%%%%%%%%%%%%%%

\title{On stable entire solutions of semi-linear elliptic equations with weights }

\author{Craig Cowan}
%    Address of record for the research reported here
\address{Department of Mathematics, Stanford University, Stanford, California.}
%    Current address
%\curraddr{Department of Mathematics and Statistics,
%Case Western Reserve University, Cleveland, Ohio 43403}
\email{ctcowan@stanford.edu}
%    \thanks will become a 1st page footnote.
\thanks{}

%    Information for second author
\author{Mostafa Fazly}
\address{Department of Mathematics, University of British Columbia, Vancouver, B.C. Canada V6T 1Z2.}
\email{fazly@math.ubc.ca}
\thanks{This work is supported by a University Graduate Fellowship and is part of the second author's Ph.D. dissertation in preparation under the supervision of N. Ghoussoub.}

%    General info

\date{..., 2011 and, in revised form, ..., 2011.}

\dedicatory{}


\date{}

\begin{abstract}  We are interested in the existence versus non-existence of non-trivial stable sub- and super-solutions of

\begin{equation} \label{pop}
-div( \omega_1 \nabla u) = \omega_2 f(u) \qquad \text{in}\ \ \IR^N,
\end{equation} with positive smooth  weights $ \omega_1(x),\omega_2(x)$. We consider the cases $ f(u) = e^u, u^p$ where $p>1$ and $ -u^{-p}$ where $ p>0$.    We obtain various non-existence results which depend on the dimension $N$ and also on $ p$ and the behaviour of $ \omega_1,\omega_2$ near infinity.  Also the monotonicity of $ \omega_1$  is involved in some results.  Our methods here are the methods developed by Farina, \cite{f2}.   We examine a specific class of weights $ \omega_1(x) = ( |x|^2 +1)^\frac{\alpha}{2}$ and $ \omega_2(x) = ( |x|^2+1)^\frac{ \beta}{2} g(x)$ where $ g(x)$ is a positive function with a finite limit at $ \infty$.   For this class of weights non-existence results are optimal.  To show the optimality we use various generalized Hardy inequalities.

\end{abstract}

\maketitle

\section{Introduction and main results}

In this note we are interested in the existence versus non-existence of stable sub- and super-solutions of equations of the form
\begin{equation} \label{eq1}
-div( \omega_1(x) \nabla u ) = \omega_2(x) f(u) \qquad \mbox{in $ \IR^N$,}
\end{equation} where $f(u)$ is one of the following non-linearities: $e^u$,  $ u^p$ where $ p>1$ and $ -u^{-p}$ where $ p>0$.  We assume that $ \omega_1(x)$ and $ \omega_2(x)$, which we call \emph{weights},  are smooth positive functions (we allow $ \omega_2$ to be zero at say a point) and which satisfy various growth conditions at $ \infty$.    Recall that we say that a solution $ u $ of $ -\Delta u = f(u)$ in $ \IR^N$ is stable provided
\[ \int f'(u) \psi^2 \le \int | \nabla \psi|^2, \qquad \forall \psi \in C_c^2,\] where $ C_c^2$ is the set of $ C^2$ functions defined on $ \IR^N$ with compact support.  Note that the stability of $u$ is just saying that the second variation at $u$ of the energy associated with the equation is non-negative.     In our setting this becomes:  We say a $C^2$ sub/super-solution $u$ of (\ref{eq1}) is \emph{stable} provided
\begin{equation} \label{stable}
\int \omega_2 f'(u) \psi^2 \le \int \omega_1 | \nabla \psi|^2 \qquad \forall \psi \in C_c^2.
\end{equation}
One should note that (\ref{eq1}) can be re-written as
\begin{equation*}
- \Delta u  + \nabla \gamma(x) \cdot \nabla u ={ \omega_2}/{\omega_1}\ f(u) \qquad \text{ in $ \mathbb{R}^N$},
\end{equation*}
where
$\gamma = - \log( \omega_1)$ and on occasion we shall take this point of view.

\begin{remark} \label{triv} Note that  if $ \omega_1$ has enough integrability then it is immediate that if $u$ is a stable solution  of (\ref{eq1}) we have $ \int \omega_2 f'(u) =0 $  (provided $f$ is increasing).   To see this let $ 0 \le \psi \le 1$ be supported in a ball of radius $2R$ centered at the origin ($B_{2R}$) with $ \psi =1$ on $ B_R$ and such that $  | \nabla \psi | \le \frac{C}{R}$ where $ C>0$ is independent of $ R$.   Putting this $ \psi$ into $ (\ref{stable})$ one obtains
\[ \int_{B_R} \omega_2 f'(u) \le \frac{C}{R^2} \int_{R < |x| <2R} \omega_1,\] and so if the right hand side goes to zero as $ R \rightarrow \infty$ we have the desired result.

\end{remark}

The existence versus non-existence of stable solutions of  $ -\Delta u =  f(u)$ in $ \IR^N$ or $ -\Delta u = g(x) f(u)$ in $ \IR^N$ is now quite well understood, see  \cite{dancer1, farina1, egg, zz, f2, f3, wei, ces, e1, e2}.  We remark that some of these results are examining the case where $ \Delta $ is replaced with $ \Delta_p$ (the $p$-Laplacian) and also in many cases the authors are interested in finite Morse index solutions or solutions which are stable outside a compact set.
  Much of the interest in these Liouville type theorems  stems from the fact that the non-existence of a stable solution is related to the existence of a priori estimates for stable solutions of a related equation on a bounded domain.

       In \cite{Ni}  equations similar to  $ -\Delta u = |x|^\alpha u^p$  where examined on the unit ball in $ \IR^N$ with zero Dirichlet boundary conditions.  There it was shown that for $ \alpha >0$ that  one can obtain positive solutions for $ p $ supercritical with respect to Sobolev embedding and so one can view that the term $ |x|^\alpha$ is restoring some compactness.      A similar feature happens for equations of the form
\[ -\Delta u = |x|^\alpha f(u) \qquad \mbox{in $ \IR^N$};\]    the value of $ \alpha$ can vastly alter the existence versus non-existence of a stable solution,  see \cite{e1, ces, e2, zz, egg}.

We now come to our main results and for this we need to define a few quantities:

\begin{eqnarray*}
I_G&:=& R^{-4t-2} \int_{ R < |x|<2R} \frac{ \omega_1^{2t+1}}{\omega_2^{2t}}dx , \\
 J_G&:=& R^{-2t-1} \int_{R < |x| <2R} \frac{| \nabla \omega_1|^{2t+1} }{\omega_2^{2t}} dx  ,\\I_L&:=& R^\frac{-2(2t+p-1)}{p-1}  \int_{R<|x|<2R }{ \left(   \frac{w_1^{p+2t-1}}{w_2^{2t}}    \right)^{\frac{1}{p-1} }  }  dx,\\ J_L&:= &R^{-\frac{p+2t-1}{p-1} }  \int_{R<|x|<2R }{ \left(   \frac{|\nabla w_1|^{p+2t-1}}{w_2^{2t}}     \right)^{\frac{1}{p-1} }  }  dx,\\
I_M  &:=&   R^{-2\frac{p+2t+1}{p+1} }  \int_{R<|x|<2R }{ \left(   \frac{w_1^{p+2t+1}}{w_2^{2t}}    \right)^{\frac{1}{p+1} }  } \ dx, \\
J_M  &:= &  R^{-\frac{p+2t+1}{p+1} }  \int_{R<|x|<2R }{ \left(   \frac{|\nabla w_1|^{p+2t+1}}{w_2^{2t}}     \right)^{\frac{1}{p+1} }  }  dx.
\end{eqnarray*}

The three equations we examine are
\[ -div( \omega_1 \nabla u ) = \omega_2 e^u \qquad \mbox{ in $ \IR^N$ } \quad  (G), \]
\[  -div( \omega_1 \nabla u ) = \omega_2 u^p \qquad \mbox{ in $ \IR^N$ } \quad  (L), \]
\[ -div( \omega_1 \nabla u ) = - \omega_2 u^{-p} \qquad \mbox{ in $ \IR^N$ } \quad  (M),\]  and where we restrict $(L)$ to the case $ p>1$ and $(M)$ to $ p>0$.    By solution we always mean a $C^2$ solution.   We now come to our main results in terms of abstract $ \omega_1 $ and $ \omega_2$.    We remark that our approach to non-existence of stable solutions is the approach due to Farina, see \cite{f2,f3,farina1}.

\begin{thm} \label{main_non_exist} \begin{enumerate}

\item   There is no  stable sub-solution of $(G)$ if $ I_G, J_G \rightarrow 0$ as $ R \rightarrow \infty$   for some   $0<t<2$.
 \item    There is no positive   stable sub-solution (super-solution) of $(L)$  if  $ I_L,J_L \rightarrow 0$ as $ R \rightarrow \infty$   for some   $p- \sqrt{p(p-1)}   < t<p+\sqrt{p(p-1)} $  ($0<t<\frac{1}{2}$).
\item   There is no positive   stable super-solution of (M) if  $ I_M,J_M \rightarrow 0$ as $ R \rightarrow \infty$  for some $0<t<p+\sqrt{p(p+1)}$.

\end{enumerate}

\end{thm}     If we assume that $ \omega_1$ has some monotonicity we can do  better.   We will assume that the monotonicity conditions is satisfied for big $x$  but really all ones needs is for it to be satisfied on a suitable sequence of annuli.

\begin{thm} \label{mono}  \begin{enumerate} \item There is no  stable sub-solution of $(G)$ with $  \nabla \omega_1(x) \cdot  x \le 0$ for big $x$ if  $ I_G \rightarrow 0$ as $ R \rightarrow \infty$   for some   $0<t<2$.

\item There is no positive  stable sub-solution of  $(L)$ provided  $  I_L \rightarrow 0$ as $ R \rightarrow \infty$ for either:
\begin{itemize}
\item  some $  1 \le t < p + \sqrt{p(p-1)}$ and $  \nabla \omega_1(x) \cdot x \le 0$ for big $x$,  or  \\

\item some $ p - \sqrt{p(p-1)} < t \le 1$ and $ \nabla \omega_1(x) \cdot x \ge 0$ for big $ x$.
\end{itemize}
 There is no positive super-solution of $(L)$ provided $ I_L \rightarrow 0$ as $ R \rightarrow \infty$ for some $ 0 < t < \frac{1}{2}$ and $ \nabla \omega_1(x) \cdot x \le 0$ for big $x$.

\item  There is no positive stable super-solution of $(M)$ provided  $ I_M \rightarrow 0$ as $ R \rightarrow \infty$  for some  $0<t<p+\sqrt{p(p+1)}$.

\end{enumerate}

\end{thm}

\begin{cor}  \label{thing}  Suppose $ \omega_1 \le C  \omega_2$  for big $ x$, $ \omega_2 \in L^\infty$,  $ \nabla \omega_1(x) \cdot  x \le 0$ for big $ x$.
\begin{enumerate} \item There is no stable sub-solution of $(G)$ if $ N \le 9$.

\item  There is no positive stable sub-solution of $(L)$ if  $$N<2+\frac{4}{p-1} \left( p+\sqrt{p(p-1)}  \right).$$

\item  There is no positive stable super-solution of $(M)$ if $$N<2+\frac{4}{p+1} \left( p+\sqrt{p(p+1)}  \right).$$

\end{enumerate}

\end{cor}

If one takes $ \omega_1=\omega_2=1$ in the above corollary, the results obtained for $(G)$ and  $(L)$,  and for some values of $p$ in $(M)$, are optimal, see \cite{f2,f3,zz}.

We now drop all monotonicity conditions on $ \omega_1$.

\begin{cor} \label{po} Suppose  $ \omega_1 \le C \omega_2$ for big $x$, $ \omega_2 \in L^\infty$, $ | \nabla \omega_1| \le C \omega_2$ for big $x$.
\begin{enumerate} \item  There is no stable sub-solution of $(G)$ if $ N \le 4$.

\item  There is no positive stable sub-solution of $(L)$ if $$N<1+\frac{2}{p-1} \left( p+\sqrt{p(p-1)}  \right).$$

\item There is no positive super-solution of $(M)$ if $$N<1+\frac{2}{p+1} \left( p+\sqrt{p(p+1)}  \right).$$

\end{enumerate}

\end{cor}

Some of the conditions on $ \omega_i$ in Corollary \ref{po} seem somewhat artificial.  If we shift over to the advection equation (and we take $ \omega_1=\omega_2$  for simplicity)
\[ -\Delta u + \nabla \gamma \cdot \nabla u = f(u), \] the conditions on $ \gamma$ become: $ \gamma$ is bounded from below and has a bounded gradient.

In what follows we examine the case where $ \omega_1(x) = (|x|^2 +1)^\frac{\alpha}{2}$ and $ \omega_2(x)= g(x) (|x|^2 +1)^\frac{\beta}{2}$,  where $ g(x) $ is positive except at say a point, smooth and where $ \lim_{|x| \rightarrow \infty} g(x) = C \in (0,\infty)$.     For this class of weights we can essentially obtain optimal results.

\begin{thm} \label{alpha_beta}   Take $ \omega_1 $ and $ \omega_2$ as above.
\begin{enumerate}

\item If $ N+ \alpha - 2 <0$ then there is no stable sub-solution for $(G)$, $(L)$ (here we require it to be positive) and in the case of $(M)$ there is no positive  stable  super-solution. This case is the trivial case, see Remark \ref{triv}.  \\

\textbf{Assumption:} For the remaining cases we assume that $ N + \alpha -2 > 0$.

  \item If  $N+\alpha-2<4(\beta-\alpha+2)$ then there is no  stable sub-solution for $ (G)$.

\item If $N+\alpha-2<\frac{ 2(\beta-\alpha+2)   }{p-1} \left( p+\sqrt{p(p-1)}  \right)$ then there is  no positive stable sub-solution of $(L)$.

\item If $N+\alpha-2<\frac{2(\beta-\alpha+2)   }{p+1} \left( p+\sqrt{p(p+1)}  \right)$ then there is no positive stable super-solution of $(M)$.

\item Further more 2,3,4 are optimal in the sense if $ N + \alpha -2 > 0$ and the remaining inequality is not satisfied (and in addition we assume we don't have equality in the inequality) then we can find a suitable function $ g(x)$ which satisfies the above properties and a stable sub/super-solution $u$ for the appropriate equation.

\end{enumerate}

\end{thm}

\begin{remark}  Many of the above results can be extended to the case of equality in either the $ N + \alpha - 2 \ge 0$ and also the other inequality which depends on the equation we are  examining.  We omit the details because one cannot prove the results in a unified way.

\end{remark}

In showing that an explicit solution is stable we will need the   weighted Hardy inequality given in \cite{craig}.
\begin{lemma} \label{Har}
Suppose $ E>0$ is a smooth function.  Then one has
\[ (\tau-\frac{1}{2})^2 \int E^{2\tau-2} | \nabla E|^2 \phi^2 + (\frac{1}{2}-\tau) \int (-\Delta E) E^{2\tau-1} \phi^2 \le \int E^{2\tau} | \nabla \phi|^2,\] for all $ \phi \in C_c^\infty(\IR^N)$ and $ \tau \in \IR$.
\end{lemma}  By picking an appropriate function $E$ this gives,

\begin{cor} \label{Hardy}
For all $ \phi \in C_c^\infty$ and $ t , \alpha \in \IR$. We have
  \begin{eqnarray*}
\int (1+|x|^2)^\frac{\alpha}{2} |\nabla\phi|^2 &\ge& (t+\frac{\alpha}{2})^2 \int |x|^2 (1+|x|^2)^{-2+\frac{\alpha}{2}}\phi^2\\
&&+(t+\frac{\alpha}{2})\int  (N-2(t+1)   \frac{|x|^2}{1+|x|^2}) (1+|x|^2)^{-1+\frac{\alpha} {2}} \phi^2.
\end{eqnarray*}
 \end{cor}

\section{Proof of main results}

\textbf{ Proof of Theorem \ref{main_non_exist}.}   (1). Suppose $ u$ is a stable sub-solution of $(G)$ with $ I_G,J_G \rightarrow 0$ as $ R \rightarrow \infty$ and  let $ 0 \le \phi \le 1$ denote a smooth compactly supported function.  Put $ \psi:= e^{tu} \phi$ into (\ref{stable}), where $ 0 <t<2$, to arrive at
\begin{eqnarray*}
\int \omega_2 e^{(2t+1)u} \phi^2 &\le & t^2 \int \omega_1 e^{2tu} | \nabla u|^2 \phi^2 \\
&& +\int \omega_1 e^{2tu}|\nabla \phi|^2 + 2 t \int \omega_1 e^{2tu} \phi \nabla u \cdot \nabla \phi.
\end{eqnarray*}  Now multiply $(G)$ by $ e^{2tu} \phi^2$ and integrate by parts to arrive at
\[ 2t \int \omega_1 e^{2tu} | \nabla u|^2 \phi^2 \le \int \omega_2 e^{(2t+1) u} \phi^2 - 2 \int \omega_1 e^{2tu} \phi \nabla u \cdot \nabla \phi,\]
and now if one equates like terms they arrive at
\begin{eqnarray} \label{start}
\frac{(2-t)}{2} \int \omega_2 e^{(2t+1) u} \phi^2  & \le & \int \omega_1 e^{2tu} \left( | \nabla \phi |^2 - \frac{ \Delta \phi}{2} \right) dx \nonumber \\
&& - \frac{1}{2} \int e^{2tu} \phi \nabla \omega_1 \cdot \nabla \phi.
\end{eqnarray}   Now substitute $ \phi^m$ into this inequality for $ \phi$ where $ m $ is a big integer to obtain
\begin{eqnarray} \label{start_1}
\frac{(2-t)}{2} \int \omega_2 e^{(2t+1) u} \phi^{2m}  & \le & C_m \int \omega_1 e^{2tu} \phi^{2m-2}  \left( | \nabla \phi |^2 + \phi |\Delta \phi|  \right) dx \nonumber \\
&& - D_m \int e^{2tu} \phi^{2m-1} \nabla \omega_1 \cdot \nabla \phi
\end{eqnarray} where $ C_m$ and $ D_m$ are positive constants just depending on $m$.   We now estimate the  terms on the right but we mention that when ones assume the appropriate monotonicity on $ \omega_1$ it is the last integral on the right which one is able to drop.

\begin{eqnarray*}
\int \omega_1 e^{2tu} \phi^{2m-2} | \nabla \phi|^2 & = & \int \omega_2^\frac{2t}{2t+1} e^{2tu} \phi^{2m-2}  \frac{ \omega_1 }{\omega_2^\frac{2t}{2t+1}} | \nabla \phi|^2  \\
& \le &  \left( \int \omega_2 e^{(2t+1) u} \phi^{(2m-2) \frac{(2t+1)}{2t}} dx \right)^\frac{2t}{2t+1}\\ &&\left( \int \frac{ \omega_1 ^{2t+1}}{\omega_2^{2t}} | \nabla \phi |^{2(2t+1)} \right)^\frac{1}{2t+1}.
\end{eqnarray*}
Now, for fixed $ 0 <t<2$ we can take $ m $ big enough so $ (2m-2) \frac{(2t+1)}{2t} \ge 2m $ and since $ 0 \le \phi \le 1$ this allows us to replace the power on $ \phi$ in the first term on the right with $2m$   and hence we obtain
 \begin{equation} \label{three}
 \int \omega_1 e^{2tu} \phi^{2m-2} | \nabla \phi|^2  \le \left( \int \omega_2 e^{(2t+1) u} \phi^{2m} dx \right)^\frac{2t}{2t+1} \left( \int \frac{ \omega_1 ^{2t+1}}{\omega_2^{2t}} | \nabla \phi |^{2(2t+1)} \right)^\frac{1}{2t+1}.
 \end{equation}    We now take the test functions $ \phi$ to be such that $ 0 \le  \phi \le 1$ with $ \phi $ supported in the ball $ B_{2R}$ with $ \phi = 1 $ on $ B_R$ and $ | \nabla \phi | \le \frac{C}{R}$ where $ C>0$ is independent of $ R$.   Putting this choice of $ \phi$ we obtain
 \begin{equation} \label{four}
 \int \omega_1 e^{2tu} \phi^{2m-2} | \nabla \phi |^2 \le \left( \int \omega_2 e^{(2t+1)u} \phi^{2m} \right)^\frac{2t}{2t+1} I_G^\frac{1}{2t+1}.
 \end{equation}  One similarly shows that
 \[ \int \omega_1 e^{2tu} \phi^{2m-1} | \Delta \phi| \le \left( \int \omega_2 e^{(2t+1)u} \phi^{2m} \right)^\frac{2t}{2t+1} I_G^\frac{1}{2t+1}.\]
 So, combining the results we obtain

 \begin{eqnarray} \label{last} \nonumber \frac{(2-t)}{2} \int \omega_2 e^{(2t+1) u} \phi^{2m} &\le& C_m \left( \int \omega_2 e^{(2t+1) u} \phi^{2m} dx \right)^\frac{2t}{2t+1} I_G^\frac{1}{2t+1}\\
 &&- D_m \int e^{2tu} \phi^{2m-1}  \nabla \omega_1 \cdot \nabla \phi.
 \end{eqnarray}
 We now estimate this last term.  A similar argument using H\"{o}lder's inequality shows that
 \[ \int e^{2tu} \phi^{2m-1} | \nabla \omega_1| | \nabla \phi| \le \left(  \int \omega_2 \phi^{2m} e^{(2t+1) u} dx \right)^\frac{2t}{2t+1} J_G^\frac{1}{2t+1}. \] Combining the results gives that
\begin{equation} \label{last}
(2-t) \left( \int \omega_2 e^{(2t+1) u} \phi^{2m} dx \right)^\frac{1}{2t+1} \le I_G^\frac{1}{2t+1} + J_G^\frac{1}{2t+1},
\end{equation} and now we send $ R \rightarrow \infty$ and use the fact that $ I_G, J_G \rightarrow 0$ as $ R \rightarrow \infty$ to see that
\[ \int \omega_2 e^{(2t+1) u} =0, \] which is clearly a contradiction.   Hence there is no stable sub-solution of $(G)$.

(2).   Suppose that $u >0$ is a stable sub-solution (super-solution) of $(L)$.  Then a similar calculation as in (1) shows that for  $ p - \sqrt{p(p-1)} <t < p + \sqrt{p(p-1)}$,  $( 0 <t<\frac{1}{2})$ one has

\begin{eqnarray}   \label{shit}
(p  -\frac{t^2}{2t-1}   )\int \omega_2 u^{2t+p-1} \phi^{2m} & \le & D_m \int \omega_1 u^{2t} \phi^{2(m-1)} (|\nabla\phi|^2  +\phi |\Delta \phi |) \nonumber \\
&& +C_m \frac{(1-t)}{2(2t-1)} \int u^{2t} \phi^{2m-1}\nabla \omega_1 \cdot  \nabla \phi.
 \end{eqnarray}    One now applies H\"{o}lder's argument as in (1) but the terms $ I_L$ and $J_L$ will appear on the right hand side of the resulting
 equation.      This shift from a sub-solution to a super-solution depending on whether $ t >\frac{1}{2}$ or $ t < \frac{1}{2}$ is a result from the sign change of $ 2t-1$ at $ t = \frac{1}{2}$.   We leave the details for the reader.

(3).  This case is also similar to (1) and (2).

\hfill $ \Box$

 \textbf{Proof of Theorem \ref{mono}.}   (1).  Again we suppose there is a stable sub-solution $u$ of $(G)$.  Our starting point  is (\ref{start_1}) and we wish to be able to drop the term
 \[ - D_m \int e^{2tu} \phi^{2m-1} \nabla \omega_1 \cdot \nabla \phi, \] from (\ref{start_1}).    We can choose $ \phi$ as in the proof of Theorem \ref{main_non_exist} but also such that $ \nabla \phi(x) = - C(x) x$ where $ C(x) \ge 0$.     So if we assume that $ \nabla \omega_1 \cdot x \le 0$ for big $x$ then we see that this last term is non-positive and hence we can drop the term.  The the proof is as before but now we only require that $ \lim_{R \rightarrow \infty} I_G=0$.

 (2).     Suppose that $ u >0$ is a stable sub-solution of $(L)$  and so (\ref{shit}) holds for all $  p - \sqrt{p(p-1)} <t< p + \sqrt{p(p-1)}$.   Now we wish to use monotonicity to drop the term from (\ref{shit}) involving the term $ \nabla \omega_1 \cdot \nabla \phi$.      $ \phi$ is chosen the same as in (1)  but here one notes that the co-efficient for this term changes sign at $ t=1$ and hence by restriction $t$ to the appropriate side of 1 (along with the above condition on $t$ and $\omega_1$) we can drop the last term depending on which monotonicity we have and hence to obtain a contraction we only require that $ \lim_{R \rightarrow \infty} I_L =0$.   The result for the non-existence of a stable super-solution is similar be here one restricts $ 0 < t < \frac{1}{2}$.

(3).  The proof here is similar to (1) and (2) and we omit the details.

 \hfill $\Box$

\textbf{Proof of Corollary \ref{thing}.}   We suppose  that $ \omega_1 \le C  \omega_2$  for big $ x$, $ \omega_2 \in L^\infty$,  $ \nabla \omega_1(x) \cdot  x \le 0$ for big $ x$.     \\
(1).    Since $ \nabla \omega_1 \cdot x \le 0$ for big $x$ we can apply Theorem \ref{mono} to show the non-existence of a stable solution to $(G)$.   Note that  with the above assumptions on $ \omega_i$ we have that
\[ I_G \le \frac{C R^N}{R^{4t+2}}.\]  For $ N \le 9$  we can take $ 0 <t<2$  but close enough to $2$ so the right hand side goes to zero as $ R \rightarrow \infty$.

Both (2) and (3) also follow directly from applying Theorem \ref{mono}.   Note that one can say more about (2) by taking the multiple cases as listed in Theorem \ref{mono} but we have choice to leave this to the reader.

\hfill $ \Box$

\textbf{Proof of Corollary \ref{po}.}   Since we have no monotonicity conditions now we will need both $I$ and $J$ to go to zero to show the non-existence of a stable solution.   Again the results are obtained immediately by applying Theorem \ref{main_non_exist}  and we prefer to omit the details.

\hfill $\Box$

\textbf{Proof of Theorem \ref{alpha_beta}.}  (1).  If $ N + \alpha -2 <0$ then using Remark \ref{triv}  one easily sees there is no stable sub-solution of $(G)$ and $(L)$ (positive for $(L)$) or a positive stable super-solution of $(M)$.   So we now assume that $ N + \alpha -2 > 0$.     Note that the monotonicity of $ \omega_1$ changes when $ \alpha $ changes sign and hence one would think that we need to consider separate cases if we hope to utilize the monotonicity results.   But a computation shows that in fact $ I$ and $J$ are just multiples of each other in all three cases so it suffices to show, say, that $ \lim_{R \rightarrow \infty} I =0$. \\
(2).   Note that for $ R >1$ one has
\begin{eqnarray*}
I_G & \le  & \frac{C}{R^{4t+2}} \int_{R <|x| < 2R} |x|^{ \alpha (2t+1) - 2t \beta} \\
& \le  &  \frac{C}{R^{4t+2}}  R^{N + \alpha (2t+1) - 2t \beta},
\end{eqnarray*} and so to show the non-existence we want to find some $ 0 <t<2$ such that
$  4t+2 > N  + \alpha(2t+1) - 2 t \beta$,   which is equivalent to  $ 2t ( \beta - \alpha +2) > (N + \alpha -2)$.    Now recall that we are assuming that $ 0 < N + \alpha -2 < 4 ( \beta - \alpha +2) $ and hence we have the desired result by taking $  t <2$ but sufficiently close.
The proof of the non-existence results for
(3) and (4) are similar and we omit the details.   \\
(5).  We now assume that $N+\alpha-2>0$.  In showing the existence of stable sub/super-solutions we need to consider  $ \beta - \alpha + 2 <0$ and $ \beta - \alpha +2 >0$ separately.

\begin{itemize} \item $(\beta - \alpha + 2 <0)$  Here we take $ u(x)=0$ in the case of $(G)$ and $ u=1$ in the case of $(L)$ and $(M)$. In addition we take $ g(x)=\E$.  It is clear that in all cases $u$ is the appropriate sub or super-solution.  The only thing one needs to check is the stability.    In all cases this reduces to trying to show that we have
\[ \sigma \int       (1+|x|^2)^{\frac{\alpha}{2}  -1}    \phi^2 \le \int      (1+|x|^2)^{\frac{\alpha}{2}}   | \nabla\phi  |^2,\]  for all $ \phi \in C_c^\infty$ where  $ \sigma $ is some small positive constant; its either $ \E$ or $ p \E$ depending on which equation were are examining.
To show this we use the result from Corollary \ref{Hardy} and we drop a few positive terms to arrive at
\begin{equation*}
\int (1+|x|^2)^\frac{\alpha}{2} |\nabla\phi|^2\ge (t+\frac{\alpha}{2})\int \left (N-2(t+1)   \frac{|x|^2}{1+|x|^2}\right) (1+|x|^2)^{-1+\frac{\alpha} {2}}
\end{equation*} which holds for all $ \phi \in C_c^\infty$ and $ t,\alpha \in \IR$.
  Now, since $N+\alpha-2>0$, we can choose $t$ such that $-\frac{\alpha}{2}<t<\frac{n-2}{2}$.  So, the integrand function in the right hand side is positive and since for small enough $\sigma$ we have
  \begin{equation*}
\sigma \le  (t+\frac{\alpha}{2})(N-2(t+1)   \frac{|x|^2}{1+|x|^2})  \ \ \ \text {for all} \ \ x\in \mathbb{R}^N
\end{equation*}
 we get stability.

\item ($\beta-\alpha+2>0$) In the case of $(G)$ we take   $u(x)=-\frac{\beta-\alpha+2}{2} \ln(1+|x|^2)$ and $g(x):= (\beta-\alpha+2)(N+(\alpha-2)\frac{|x|^2}{1+|x|^2})$. By a computation one sees that $u$ is a sub-solution of $(G)$ and hence we need now to only show the stability, which amounts to showing that
\begin{equation*}
\int \frac{g(x)\psi^2}{(1+|x|^{2   })^{-\frac{\alpha}{2}+1}}\le \int\frac{|\nabla\psi|^2}{    (1+|x|^2)^{-\frac{\alpha}{2}}     },
\end{equation*} for all $ \psi \in C_c^\infty$.  To show this we use  Corollary \ref{Hardy}.  So  we  need to choose an appropriate $t$ in   $-\frac{\alpha}{2}\le t\le\frac{N-2}{2}$  such that for all $x\in \IR^N$ we have
 \begin{eqnarray*}
 (\beta-\alpha+2)\left(    N+  (\alpha-2)\frac{|x|^2}{1+|x|^2}\right)         &\le& (t+\frac{\alpha}{2})^2 \frac{ |x|^2 }{(1+|x|^2}\\
&&+(t+\frac{\alpha}{2}) \left(N-2(t+1)   \frac{|x|^2}{1+|x|^2}\right).
\end{eqnarray*}
With a  simple calculation one sees we need just to have
   \begin{eqnarray*}
 (\beta-\alpha+2)&\le& (t+\frac{\alpha}{2}) \\
  (\beta-\alpha+2)     \left(    N+  \alpha-2\right)      &   \le&  (t+\frac{\alpha}{2}) \left(N-t-2+\frac{\alpha}{2}) \right).
 \end{eqnarray*}     If one takes $ t= \frac{N-2}{2}$ in the case where $ N \neq 2$ and $ t $ close to zero in the case for $ N=2$ one easily sees the above inequalities both hold, after considering all the constraints on $ \alpha,\beta$ and $N$.

 We now consider the case of $(L)$.  Here one takes $g(x):=\frac {\beta-\alpha+2}{p-1}(    N+  (\alpha-2-\frac{\beta-\alpha+2}{p-1})
\frac{|x|^2}{1+|x|^2})$ and $ u(x)=(1+|x|^2)^{ -\frac  {\beta-\alpha+2}{2(p-1)} }$.   Using essentially the same approach as in $(G)$ one shows that $u$ is a stable sub-solution of $(L)$ with this choice of $g$.   \\
For the case of $(M)$ we take   $u(x)=(1+|x|^2)^{ \frac  {\beta-\alpha+2}{2(p+1)}   }$ and $g(x):=\frac {\beta-\alpha+2}{p+1}(    N+  (\alpha-2+\frac{\beta-\alpha+2}{p+1})
\frac{|x|^2}{1+|x|^2})$.

\end{itemize}

\hfill  $  \Box$

 \end{document}